\documentclass[12pt]{article}
\usepackage{amsmath,amsfonts,amssymb}
\usepackage{graphicx,color}
\usepackage{theorem}

\theorembodyfont{\upshape}

\newcommand{\qed}{\nopagebreak\par\hspace*{\fill}$\square$}

\newcommand{\ZZ}{{\mathbb Z}}

\newcommand{\CC}{{\mathbb C}}

\def\gg{{\mathfrak{g}}}

\def\hh{{\mathfrak{h}}}

\def\pp{{\mathfrak{p}}}
\def\bb{{\mathfrak{b}}}
\def\ss{{\mathfrak{s}}}
\def\aa{{\mathfrak{a}}}

\numberwithin{equation}{section}

\title{A new class of modules for Toroidal Lie Superalgebras}
\author{S.Eswara Rao\\
School of Mathematics\\
Tata Institute of Fundamental Research,\\
Mumbai, India.\\[2mm]
email: senapati@math.tifr.res.in}
\date{}

\begin{document}

\maketitle

\begin{abstract}
In this paper we construct a large class of modules for toroidal Lie superalgebras. 
Toroidal Lie superalgebras are universal central extensions of $\gg \otimes A$
where $\gg$ is a basic classical Lie superalgebra and $A$ is Laurent polynomial
ring in several variables. The case where $\gg$ is a simple finite dimensional 
Lie algebra is included.
\end{abstract}

\section*{Introduction}
The purpose of this paper is to construct a large class of modules for toroidal
Lie superalgebra. Toroidal Lie superalgebras are universal central extensions of 
$\gg \otimes A$ where $\gg$ is a basic classical simple Lie superalgebra and $A$ is a 
Laurent polynomial ring in several variables. These algebras are first studied in $[IK]$ and
$[EZ]$. When $\gg$ is a simple finite dimensional Lie algebra we get toroidal Lie 
algebras which are extensively studied. Frenkel, Jing and Wang $[FJW]$ use representations
of toroidal Lie algebras to construct a new form of the McKay correspondence. Inami et al.
studied toroidal symmetry in the context of 4-dimensional conformal theory $[IKUX] , [IKU]$.
There are also applications of toroidal Lie algebras to soliton theory. Using representations
of toroidal Lie algebras one can construct hierchies of non-linear $PDEs [B], [ISW].$
We hope similar application can be found for Toroidal Lie superalgebras.

In this paper we construct a functor from modules of affine superalgebras to modules of toroidal Lie superalgebras.
Even in the Lie algebra case our construction is completely new and produce a large class of modules
for toroidal Lie algebras. Our construction recover the well known results of $[EMY]$ and $[EM].$ 
In $[EM]$ and $[EMY]$, only the level one integrable modules are considered but our construction
works for any highest weight modules and any non-zero level.

We will now explain the results in more detail. Let $\gg$ be a basic classical Lie superalgebra and 
$A = \CC [t_1^{\pm 1},\cdots, t_n^{\pm1}]$ be a Laurent  polynomial ring in $n$ commuteing 
variables. Than $\gg \otimes A$ is naturally a Lie superalgebra. The universal central extension 
$\tau$ of $\gg \otimes A$ is called toroidal Lie superalgebra. It is explicitly given in $[IK]$. 
Let $\gg_{aff}$ be the affine superalgebra corresponding to $\gg$ (See 2.4). We first construct the 
standard Fock space $V(\Gamma)$ for a suitable non-degenerate lattice $\Gamma$ and a degenerate
sublattice $Q$ (see section 3). Let $V$ be any restricted $\gg_{aff}$ module (See 2.5) and consider $V \otimes V(\Gamma)$. 
We prove the following
\begin{enumerate}
\item $V \otimes V(\Gamma)$ is a $\stackrel{\sim}{\tau}$- module (Theorem 4.1)
\item Suppose $V$ is a weight module for $\gg_{aff}$ then $V \otimes V(\Gamma)$ is a weight module for 
 $\stackrel{\sim}{\tau}$ for suitable Cartan subalgebra.
\item Suppose $V$ is a $\gg_{aff}$ module in category ${\cal O}$ with finite dimensional weight spaces then $V \otimes V(\Gamma)$
is a weight module with finite dimensional weight spaces (Remark 6.2).
\item Suppose $\gg$ is a Lie algebra and suppose $V$ is a  $\gg_{aff}$ integrable module then $V \otimes V(\Gamma)$ is
$\stackrel{\sim}{\tau}$ - integrable.
\item Suppose $V$ and $W$ are restricted $\gg_{aff}$ modules and $ f : V \rightarrow W$ is $\gg_{aff}$ module map.
Then there exists a canonical map $\stackrel{\sim}{f} : V \otimes V(\Gamma) \rightarrow W \otimes V(\Gamma)$ which is a 
$\stackrel{\sim}{\tau}$ - module map (Remark 6.3).
\end{enumerate}
We further note that $V= \displaystyle{\bigoplus_{\lambda \in \Gamma/Q}} V (\lambda)$ where each $V (\lambda)$ is a
$\stackrel{\sim}{\tau}$ - module. We prove in Proportion $5.2$ that $V 
(\lambda)$ and $V (0)$ are isomorphic as
$\stackrel{\sim}{\tau}$ - modules upto a twist of an automorphism of 
$\stackrel{\sim}{\tau}$.

\section {Lie Superalgebra} 
A Lie superalgebra is a $\ZZ_{2}$- graded vector space $\gg = \gg_{\overline{0}} \oplus \gg_{\overline{1}}$ equipped with $\CC$- bilinear form
$[,] : \gg \times \gg \rightarrow \gg$, called the Lie super bracket, satisfying the following conditions.
\begin{enumerate}
\item  $[\gg_{\overline{i}}, \gg_{\overline{j}}] \ {\underline{\subset}} \  \gg_{\overline{i+j}}$

\item $[X,Y] =  - (-1)^{ij} \ [Y,X]$

\item $\big[[X,Y],Z\big] = \big[X,[Y,Z]\big] - (-1)^{ij} \ \big[Y,[X,Z]\big]$
\end{enumerate}
for all homogeneous elements $X \in \gg_{\overline{i}}, Y \in \gg_{\overline{j}}$ and $\ZZ \in \gg_{\overline{k}}$.

The subspace $\gg_{\overline{0}}$ is called even and the subspace $\gg_{\overline{1}}$ is called odd. It is easy to 
see that $\gg_{\overline{0}}$ is the usual Lie algebra and $\gg_{\overline{1}}$ is $\gg_{\overline{0}}$ - module. The 
identity $(3)$ is called super Jacobi identity. Suppose $X$ is a homogeneous element belonging to $\gg_{\overline{i}}$, 
then we denote $|X| = i$.

A bilinear form $(,) : \gg \times \gg \rightarrow \CC$ is called
\begin{enumerate}
\item  Supersymmetric if $(X,Y) = (-1)^{|X| |Y|} (Y,X)$ for all homogeneous elements $X$ and $Y$ in $\gg$.

\item Invariant if $\big([X,Y], Z\big) = \big(X,[Y,Z]\big)$ for all $X,Y,Z \in \gg$.

\item Even if $(X,Y) = 0$ for all $ X \in \gg_{\overline{0}}$ and $Y \in \gg_{\overline{1}}$.
\end{enumerate}

A Lie superalgebra $\gg$ is called basic classical if $\gg$ is simple, finite dimensional, the even part is 
reductive and $\gg$ carries an even, non-degenerate supersymmetric bilinear form. They have been classified by Kac
$[K2]$. The following is the list of basic classical Lie superalgebra and the decomposition of the even part.
$$
\begin{array}{llll}
A(m,n)& A_{m} + A_{n} + \CC& ,m \geq 0,n \geq 0, m+n \geq 1\\
B(m,n)& B_{m} + C_{n}& ,m \geq 0,n \geq 1\\
C(n)& C_{n-1} + \CC& ,n \geq 3\\
D(m,n)& D_{m} + C_{n}& ,m \geq 2,n \geq 1\\
D(2,1,a)& D_{2} + A_{1}& a\neq 0, -1\\
F(4)& B_{3}+A_{1}\\
G(3)& G_{2}+A_{1}\\
\end{array}
$$
In addition to this we need to add all simple finite dimensional Lie algebras.

Suppose $\gg$ is a basic classical Lie superalgebra and $(,)$ be an even, non-degenerate, supersymmetric and invariant
bilinear form. Let $\hh \ \underline{\subset} \ \gg$ be a Cartan subalgebra.

Let $0 \neq \alpha \in \hh^{*}$ and\\ let
$\gg_{\alpha} = \{X \in \gg \ | \ [h,X] = \alpha (h) X, \forall h \in \hh \}$,
Let $\Delta = \{\alpha \in \hh^{*} \ | \ \gg_{\alpha} \neq 0 \}$\\
Then it is well known that
$$
\gg = \displaystyle {\bigoplus_{\alpha \in \Delta}} \ \gg_{\alpha} \oplus \hh.
$$

The form $(,)$ restricted to $\hh$ is non-degenerate. Using this form, $\hh^{*}$ can be identified with $\hh$ 
via the map $\alpha \mapsto H (\alpha)$\\
Let $\Delta_{0} = \{ \alpha \in \Delta \ |\ \gg_{\alpha} \cap \gg_{\overline{0}} \neq 0 \}$ and\\
$\Delta_{1} = \{ \alpha \in \Delta \ |\ \gg_{\alpha} \cap \gg_{\overline{1}} \neq 0 \}$.\

Let $\Delta^{+}$ and $\Delta^{-}$ denote the positive and negative roots respectively.\\
Let $\Delta_{1}^{+} = \Delta_{1} \cap \Delta^{+}, \Delta_{1}^{-}= \Delta_{1} \cap \Delta^{-}.$\\
For $\alpha \in \Delta_{1}^{-}$, define $\sigma_{\alpha} = -1$ and $\sigma_{\alpha} = 1$ for $\alpha \notin \Delta_{1}^{-}$.\\
Let $\dim \hh = l$.
\paragraph{1.1}
Let $\gg$ be basic classical Lie superalgebra. We call a Chevelley basis $\gg$, any homogeneous $\CC$ -basis of $\gg$,\\
$B= \{ H(i)\}_{i=1,2,\cdots l} \ \cup \{X(\alpha), \alpha \in \Delta \}$\\
such that
\begin{enumerate}
\item[{(a)}] $H(1), \cdots H(l)$ is a $\CC$ basis of $\hh$

\item[{(b)}] $[H(i), H(j)] =0 , [H(i), X(\alpha)] = \alpha (H(i)) X (\alpha)$ for $i,j=1,2,\cdots l$ and $\alpha \in \Delta$.

\item[{(c)}] $[X(\alpha), X(-\alpha)] = \sigma_{\alpha} H(\alpha), \alpha \in \Delta$

\item[{(d)}] $[X(\alpha), X(\beta)] = N_{\alpha, \beta} \ X(\alpha + \beta) , \forall \ \alpha, \beta \in \Delta$

\item[{(d1)}] $N_{\alpha, \beta} \in \ZZ$ and $N_{\alpha, \beta} = 0$ for $\alpha + \beta \notin \Delta \cup \{0\}$

\item[{(d2)}] If $(\alpha, \alpha) \neq 0$\ or $(\beta,\beta) \neq 0$ and if \\
$\sum^{\alpha}_{\beta} = \{ \beta - r \alpha, \cdots \beta + q \alpha \}$ is the $\alpha$ string through $\beta$ 
then\\ $N_{\alpha, \beta} = \pm ( r + 1).$

\item[{(d3)}] If $(\alpha, \alpha) = 0 = (\beta,\beta)$ then\\ $N_{\alpha, \beta} = \pm \ \beta \ (H(\alpha))$
\end{enumerate}

\paragraph{(1.2) Proposition} $\big([IK],[FG]\big).$ Every basic classical Lie superalgebra admits a Chevelley basis.

\paragraph{(1.3) Remark} Suppose $\gg = sl (d+1,d+1)$ which is not a basic classical Lie superalgebra. But it is known
that Chevelley basis exists. See $[IK]$.

\section{Toroidal Superalgebra}
In this section we define toroidal superalgebra and fix some notation. We fix a positive integer $n$. Let 
$A=A_{n} = \CC [t_1^{\pm 1},\cdots, t_n^{\pm1}]$ be a Laurent polynomial ring in $n$ commuteing variables.
Let ${\overline{m}} = (m_{1}, \cdots, m_{n}) \in \ZZ^{n}$ and let \\
$t^{\overline{m}} = t^{{m_1}}_{1} \cdots t^{m_{n}}_{n} \in A.$ Let $\gg$ be a basic classical Lie superalgebra
and we fix an even, non-degenerate supersymmetric invariant bilinear form $(,)$on $\gg$. For any vector space $V$ 
over $\CC$ we denote $V \otimes A$ by $V_{A}$ and $v ({\overline{m}}) = v \otimes t^{\overline{m}} \in V_{A}$

Then $\gg \otimes A$ has a natural Lie superalgebra structure. Let $\Omega_{A}/ d_{A}$ be a space of differentials so that 
$\Omega_A$ is spanned by symbols $t^{\overline{m}} K_{i}, 1 \leq i \leq n, {\overline{m}} \in \ZZ^{n}$ and $d_A$ is spanned
by $\sum  m_{i} t^{\overline{m}} K_{i}.$ We define Lie superalgebra structure on 
$\tau = \gg \otimes A \oplus \ \Omega_{A}/ d_{A}$

\paragraph{(2.1)} 
$[X({\overline{m}}), Y({\overline{k}})] = [X,Y] ({\overline{m}}+{\overline{k}})+(X,Y) d  (t^{\overline{m}})  t^{\overline{k}}$\\
where $d  (t^{\overline{m}})  t^{\overline{k}} = \displaystyle{\sum_{i}} m_{i} t^{{\overline{m}}+{\overline{k}}} K_{i}, X,Y \in \gg, {\overline{m}},
{\overline{k}} \in \ZZ^{n} ,\ \Omega_{A}/ d_{A}$ is central in $\tau$.
\paragraph{(2.2) Remark.} The above construction holds good for the Lie superalgebra $sl(d+1,d+1), d \geq1.$
The canomical form is degenerate and has a one-dimensional radical.
\paragraph{(2.3) Theorem} (Theorems (4.7) of [IK]) $\tau$ is the universal central extension of $\gg \otimes A.$
(In the case $\gg$ is of type $A(d,d)$ we take $\gg$ to be $sl(d+1,d+1)).$\\

$\tau$ is naturally $\ZZ^{n}$-graded and to 
reflect this fact we add a finite set of derivations. Let $D$ be the vector space spanned by $d_{1}, \cdots d_{n}$ and let\\
$\stackrel{\sim}{\tau} = \tau \oplus D.$\\
Define\\
$[d_{i}, X({\overline{m}})] = m_{i}  X({\overline{m}}), X \in \gg, {\overline{m}} \in \ZZ^{n}$\\
$[d_{i},d  (t^{\overline{m}})  t^{\overline{k}}]= (m_{i}+k_{i}) d  (t^{\overline{m}})  t^{\overline{k}}$\\
$[d_{i}, d_{j}] = 0$\\
Then $\stackrel{\sim}{\tau}$ becomes a Lie superalgebra with even part  
$\gg_{\overline{0}} \otimes A \oplus \Omega_{A}/ d_{A} \oplus D$ and odd part of $\gg_{\overline{1}} \otimes A.$ 
Let $\hh$ be a Cartan subalgebra of $\gg$ which is contained in $\gg_{\overline{0}}$. 
Then $\stackrel{\sim}{\hh} = \hh \oplus \displaystyle{\sum_{i=1}^{n}} \ \CC  K_{i} \oplus D$ is a Cartan subalgebra of $\stackrel{\sim}{\tau}$.\\
For $1\leq i \leq n,$ let $\delta_{i} \in \stackrel{\sim}{\hh}^{*}$ defined by $\delta_{i}(\hh)=0, \delta_{i}(K_{j}) =0$ and
$\delta_{i}(d_{j}) = \delta_{ij}.$ Let $\delta_{\overline{m}} = \displaystyle{\sum_{i}} m_{i} \delta_{i}$.\\
Let $\tau_{{\alpha} + \delta_{\overline{m}}}= \gg_{\alpha} \otimes t^{\overline{m}}, \alpha \in \Delta$,
$\tau_{\delta_{\overline{m}}} = \hh \otimes t^{\overline{m}} , {\overline{m}} \neq 0$,
$\tau_{0} = \stackrel{\sim}{\hh}$.\\
Then $\stackrel{\sim}{\tau} = \displaystyle {\bigoplus_{\substack{\alpha \in \Delta\\{\overline{m}} \in \ZZ^{n}}}} \tau_{\alpha +\delta_{\overline{m}}}$ 
is a root space decomposition with respect to the subalgebra $\stackrel{\sim}{\hh}$. See $[EZ]$ for more details.

\paragraph{(2.4) Affine superalgebra} 
Let $\gg_{aff}= \gg \otimes \CC [t_{n},t^{-1}_{n}] \oplus \CC K_{n} \oplus \CC {\overline{d}}_{n}$\\ and the Lie bracket is given by
$[X \otimes t^{m}_{n}, Y \otimes t^{k}_{n}] = [X,Y] \otimes t^{m+k}_{n} + (X,Y) m \ \delta_{m+k,0} K_{n}$.\\
$K_n$ is central and $[{\overline{d}}_{n}, X \otimes t^{m}_{n}] = m X \otimes t^{m}_{n}.$\\
Then $\gg_{aff}$ is called Affine superalgebra corresponding to the basic classical superalgebra.\

Let $\delta_{n}$ be the null root. Then $\Delta_{aff} = \{\alpha + m \delta_{n}, m \delta_{n}, \alpha \in \Delta, m \in \ZZ\}$ 
is a root system for $\gg_{aff}$.

\paragraph{(2.5)} A $\gg_{aff}$ module $V$ is called restricted if for any $v$ in $V$,
$X \otimes  t^{m}_{n} v = 0$ for $m > > 0, X \in \gg$.\

The purpose of this paper is to construct a functor from restricted $\gg_{aff}$ modules to $\stackrel{\sim}{\tau}$ modules.
The functor takes weight modules to weight modules and integrable modules to integrable modules (In the Lie algebra case. See Remark (6.1)).
But an irreducible modules never goes to irreducible module. It goes to an indecomposable module. This constructions works for any
finite dimensional simple Lie algebra and a completely new result. This construction recovers the well known result in the 
paper $[EM]$ and $[EMY].$

\section{Fock Space}
\paragraph{(3.1)}Let $\Gamma$ be free $\ZZ-$module on genarators $\delta_{1},\delta_{2},\cdots \delta_{n-1}$ and $d_{1},d_{2},\cdots d_{n-1}$
and define a symmetric bilinear form on $\Gamma$ by $< d_{i}, d_{j}> = < \delta_{i},\delta_{j} > = 0$ and 
 $< \delta_{i}, d_{j} > = \delta_{ij}$ so that 
$$
\Gamma = \displaystyle{\bigoplus_{i=1}^{n-1}} \ \ZZ  \delta_{i} \ {\bigoplus_{i=1}^{n-1}} \ \ZZ d_{i}
$$
Let $Q= {\oplus_{i=1}^{n-1}} \ \ZZ \delta_{i}$,
Let $\pp = \CC \otimes_{\ZZ} \Gamma$ and $\ss = \CC \otimes_{\ZZ} Q$,
$\bb = \oplus_{k \in \ZZ} \ \pp (k) \oplus \CC  c$\\
Where each $\pp(k)$ is an isomorphic copy of $\pp$ and the isomorphism is given by $\beta \mapsto \beta(k).$
The Lie algebra structure is given by 
$$
[\alpha (k), \beta(m)] = k <\alpha,\beta> \delta_{m+k,0} \ c
$$
and $c$ is central.\\
Let $\bb_{+}= \oplus_{k>0} \ \pp (k)$ and $\bb_{-}= \oplus_{k<0} \ \pp (k)$,
so that $\bb= \bb_{+} \oplus \bb_{-} \oplus \CC \ c$\\
Similary define $\aa ,\aa_{\pm}$ by replacing $\pp$ by $\ss$\

The Fock space representation of $\bb$ is the symmetric algebra $S(\bb_{-})$ of $\bb_{-}$ together with the action of 
 $\bb_{-}$ on  $S(\bb_{-})$ defined by\

$c$ acts on $Id.$\

$a(-m)$ acts as multiplication by $a(-m), m>0$.\

$a(m)$ acts as unique derivations on $S(\bb_{-})$ for which \\$b(-k) \rightarrow \delta_{m,k}\  m <a, b>, m >0, k>0.$
For each $\gamma$ in $\Gamma$ let $e^\gamma$ be a symbol and form the vector space \\
$\CC [\Gamma] =\displaystyle{\sum_{\gamma \in \Gamma}} \CC e^{\gamma}$ over $\CC$. In particuler
$\CC [\Gamma]$ contains the subspace \\$\CC [Q] = \displaystyle{\sum_{\alpha \in Q}} \CC e^\alpha$.
Define multiplication on $\CC [\Gamma]$ by $e^{\alpha}. e^{\gamma} = e^{\alpha + \gamma}, \alpha,\gamma \in \Gamma.$\\
Let $M \subseteq S(\bb_-)$ be any $\aa$ submodule (with respect to Fock space action.) and let $V(\Gamma,M) = \CC [\Gamma] \otimes M.$\

Of particuler interest in the sequal will be $V(\Gamma, S(\aa_{-}))$ and $V(\Gamma, S(\bb_{-}))$ which we will simply 
denote by $V(\Gamma)$ and $V(\Gamma, \bb).$ We extend the action of $\aa$ on $M$ to $\stackrel{\wedge}{\aa}$ on $V(\Gamma, M)$ by
$$
\begin{array}{llll}
a(m). e^{\gamma} \otimes u &=& e^{\gamma} \otimes a (m)u, m \neq 0\\
a(0). e^{\gamma} \otimes u &=& (\gamma,a) e^{\gamma} \otimes u.
\end{array}
$$

\paragraph{(3.2) Vertex Operators.} Let $z$ be a complex variable  and let $\alpha \in Q$\\
Define $T_{\pm} (\alpha,z) = -\displaystyle{\sum_{n\gtrless0}} \ \frac{1}{n} \ \alpha(n) \ z^{-n}$.\\
Then the vertex operator $X(\alpha,z) = exp \ T (\alpha,z)$ where \\
$exp \ T (\alpha,z) = exp \ T_{-}(\alpha,z) e^{\alpha} \ z^{\alpha (0)} exp \ T_{+}(\alpha,z)$.\\
The operator $z^{\alpha(0)}, e^{\alpha}$ defined as 
$$
\begin{array}{llll}
z^{\alpha(0)}. e^{\gamma} \otimes u &=& z^{(\alpha,\gamma)} e^{\gamma} \otimes u\\
e^{\alpha}. e^{\gamma} \otimes u &=& e^{\alpha+\gamma} \otimes u.
\end{array}
$$
Write $X(\alpha,z) = \displaystyle{\sum_{m \in \ZZ}} X_{m} (\alpha) z^{-m}$\\
It is standard fact that $X_{m} (\alpha)$ act on $F$ and for any $v$ in $F ,X_{m} (\alpha) v =0$ for $m>>0$.\\
It is easy to see that
\paragraph{(3.3)} $X(\alpha,z) X(\beta,z)=X(\alpha+\beta,z)$ for $\alpha,\beta \in Q$ as $(\alpha,\beta)=0$.\\
The following Lemma is very standard. See For example $[FK]$ or $[EM]$

\paragraph{(3.4) Lemma}
$$
\begin{array}{llll}
(1)& [b(k),X_{m} (\alpha)] = <b,\alpha > X_{m+k} (\alpha),b \in \pp, \alpha \in Q \ ,m,k \in \ZZ\\
(2)& [X_{m} (\alpha), X_{k} (\beta)] = 0, \alpha,\beta \in Q\  m,k \in \ZZ
\end{array}
$$
Define for $h \in \pp$
$$
\begin{array}{llll}
h(z) = \displaystyle{\sum_{m \in \ZZ}} \ h (m) z^{-m-1}\\
h^{+}(z) = \displaystyle{\sum_{m \geq 0}} \ h (m) z^{-m-1}\\
h^{-}(z) = \displaystyle{\sum_{m < 0}} \ h (m) z^{-m-1}\\
\end{array}
$$
So that $h(z) = h^{+}(z) + h^{-}(z)$.\\
For $\alpha \in Q, h \in \pp$ define \\
$T^{h} (\alpha,z) = h^{-} (z) X (\alpha,z) + X(\alpha,z) \ h^{+} (z)$\\
and write $T^{h} (\alpha,z) = \displaystyle{\sum_{k \in \ZZ}} \ T^{h}_{k} (\alpha) z^{-k-1}$
\paragraph{(3.5)} Recall that $<,>$ is a non-degenerate form on $\pp$. Let $\{\alpha_{i}\}, \{\alpha^{i}\}$ be a 
dual basis for $\pp$ so that $<\alpha_{i}, \alpha^{j}> = \delta_{ij}$\\
Define 
$$
L_{0}= -\frac{1}{2} \ \sum_{i} \ \sum_{k \in \ZZ} : \alpha_{i} (k) \alpha^{i} (-k) :
$$
Where the normal ordering is defined on $:\alpha(-k) \ \beta(k):$\\
$= \alpha(-k) \beta(k)$ if $k \geq 0$\\
$=\beta(k) \alpha(-k)$ if $k < 0.$\\
Then the following is very standard. See $[EM].$
$$
\begin{array}{llll}
(1)L_{0}. e^{\gamma} \otimes a_{1}(-k_{1}) \cdots a_{m} (-k_{m})= -(\frac{(\gamma,\gamma)}{2} + k_{1} + \cdots + k_{m}).\\
e^{\gamma} \otimes a_{1} (-k_{1}) \cdots a_{m} (-k_{m})
(\gamma \in \Gamma, a_{i} \in  \pp, k_{i} \in \ZZ)\\
(2)\left[L_{0}, X_{k} (\delta_{\underline{m}})\right] = kX_{k} (\delta_{\underline{m}})
\end{array}
$$
Let ${\underline{m}}=(m_{1},\cdots m_{n-1}) \in \ZZ^{n-1}, {\overline{m}} = ({\underline{m}}, m_{n}), \ 
\delta_{\underline{m}} = \displaystyle{\sum^{n-1}_{i=1} \ m_{i}} \ \delta_{i}$
\paragraph{(3.6) Lemma} (Lemma (3.13) of $[EM]$)\\
$T^{\delta_{\underline{m}}}_{m_{n}} \ (\delta_{\underline{m}})+ m_{n} \ X_{m_{n}} (\delta_{\underline{m}}) = 0$\\
An easy way to see this is to differentiate $X(\delta_{\underline{m}},z)$ and compare coefficents.\\

We will now introduce delta function and recall some standard facts from section 2 of $[FLM].$\\
Define delta function
$$
\delta (z) = \sum_{n \in \ZZ} \ z^{n}.
$$
(This delta function is not be confused with $\delta_{\underline{m}} (z)$ as the later always comes with index).\\
Then the following Lemma holds. See $[FLM]$ for proof and definitions.\\
Suppose $X(z,w)= \displaystyle{\sum_{m,k \in \ZZ}} \ X_{m,k} \ z^{m} w^{k},$\\
define $D_{z} X(z,w)= \displaystyle{\sum_{m,k \in \ZZ}}m X_{m,k} \ z^{m-1} w^{k}$
\paragraph{(3.7) Lemma}
$$
\begin{array}{llll}
(1)& X(z,w) \delta (z/w) &=& X(w,w) \delta(z/w)\\
(2)& X(z,w)  D_{z} (\delta (z/w)) &=& X(w,w)  D_{z} (\delta (z/w)) - (D_{z} X) (w,w)  \delta (z/w)
\end{array}
$$
We need the following
\paragraph{(3.8) Lemma}
$$
\begin{array}{llll}
(1)& D_{z}  z  X (\delta_{\underline{m}},z) &=& z \ \delta_{\underline{m}}(z)  X (\delta_{\underline{m}},z) + X (\delta_{\underline{m}},z)\\
(2)& D_{z}  X (\delta_{\underline{m}},z) &=& \delta_{\underline{m}} (z)  X (\delta_{\underline{m}},z)
\end{array}
$$
{\bf{Proof}} Follows from Lemma (4.6) of $[EM]$. Just note that $\alpha(z)$ in our paper and $\alpha(z)$ in $[EM]$ differ by a
$z$ factor.\\
We need the following Lemma for later use.
\paragraph{(3.9) Lemma}
$$
\begin{array}{llll}
(1)& D_{z}  \delta(z/w). z\ X (\delta_{\underline{m}},z)  X (\delta_{\underline{k}},w)\\
&= w X (\delta_{\underline{m}+{\underline{k}}},w) D_{z}  \delta(z/w)\\
&- w  \delta_{\underline{m}} (w) X (\delta_{\underline{m}+{\underline{k}}},w)  \delta(z/w)\\
&- X (\delta_{\underline{m}+{\underline{k}}},w)  \delta(z/w)\\[2mm]
(2)& D_{z}  \delta(z/w) X (\delta_{\underline{m}},z)  X (\delta_{\underline{k}},w)\\
&=  X (\delta_{\underline{m}+{\underline{k}}},w) D_{z}  \delta(z/w)\\
&- \delta_{\underline{m}} (w) X (\delta_{\underline{m}+{\underline{k}}},w)  \delta(z/w)\\
\end{array}
$$
\section{Representations of Toroidal Superalgebra}
Let $\gg$ be a basic classical Lie superalgebra. $\big($When $\gg$ is of type $A (d,d),$ we will allow $\gg$ to be 
$sl(d+1, d+1)\big).$ Let $\gg_{aff}$ be the corresponding super affine Lie-algebra. Recall that $V(\Gamma,M)$ is a 
$\stackrel{\wedge}{\aa}$-module where $X(\delta_{\underline{m}},z)$ and $T^{\delta_{\underline{k}}} (\delta_{\underline{m}},z)$
act.\

We will fix a Chevelley basis $\{H(i)\}_{i=1,\cdots l} \ \cup\{X(\alpha)\}_{\alpha \in \Delta}$ for $\gg$.\\
Let $X_{k} (\alpha) = X (\alpha) \otimes t^{k}$ and $H_{k} (i) = H(i) \otimes t^{k}$.\\
Let $V$ be a restricted $\gg_{aff}$ module where the canonical central element $K_{n}$ act as a non-zero scalar $K$.
\paragraph{(4.1) Theorem} Notation as above.\\
Then  $V \otimes V(\Gamma)$ is a $\stackrel{\sim}{\tau}$-module under the following map.
$$
\begin{array}{llll}
X (\alpha) \otimes t^{\overline{m}}& \longmapsto & \sum X_{k} (\alpha) \otimes X_{m_{n}-k} (\delta_{\underline{m}})\\
H (i) \otimes t^{\overline{m}}& \longmapsto & \sum H_{k} (i) \otimes X_{m_{n}-k} (\delta_{\underline{m}})\\
\end{array}
$$
$1 \leq i \leq n-1$
$$
\begin{array}{llll}
t^{m} K_{i} &\longmapsto& K \otimes T^{\delta_{i}}_{m_{n}} (\delta_{\underline{m}})\\
t^{m} K_{n} &\longmapsto& K \otimes X_{m_{n}} (\delta_{\underline{m}})\\
\end{array}
$$
$1 \leq i \leq n-1$
$$
\begin{array}{llll}
d_{i} &\longmapsto& 1 \otimes d_{i} (0)\\
d_{n} &\longmapsto& {\overline{d}_{n}} \otimes 1 + 1 \otimes L_{0}\\
\end{array}
$$
\noindent
{\bf{Notation}} $\displaystyle{\sum_{k}}$ we mean the summation over all integers.\\[3mm]
{\bf{Proof}} In terms of infinite series we can write the map in the following way 
$$
\begin{array}{llll}
\displaystyle{\sum_{k}} X (\alpha) \otimes t^{\underline{m}} \ t^{k}_{n} \ z^{-k}& \longmapsto &  X (\alpha,z) X(\delta_{\underline{m}},z)\\
\displaystyle{\sum_{k}} H (i) \otimes t^{\underline{m}} \ t^{k}_{n} \ z^{-k-1}& \longmapsto &  H (i,z) X(\delta_{\underline{m}},z)\\
\end{array}
$$
$1 \leq i \leq n-1$
$$
\begin{array}{llll}
\displaystyle{\sum_{k}}  t^{\underline{m}} \ t^{k}_{n} \ K_{i} \  z^{-k-1}& \longmapsto & K \delta_{i} (z) X (\delta_{\underline{m}},z)\\
\displaystyle{\sum_{k}}  t^{\underline{m}} \ t^{k}_{n} \ K_{n} \  z^{-k}& \longmapsto & K  X (\delta_{\underline{m}},z)\\
\end{array}
$$
Where $X(\alpha,z) = \displaystyle{\sum_{k}} \ X_{k} (\alpha) z^{-k},$\\
$H(i,z) = \displaystyle{\sum_{k}} \ H_{k} (i) z^{-k-1}$
are operators acting on the affine module $V$.\\
In order to prove the Theorem, it is sufficiant to verify the following relations as operators acting on 
$V \otimes V(\Gamma,M)$. When we write these infinite series in components, we will see that they satisfy 
bracket operations of $\stackrel{\sim}{\tau}$.\\

\noindent
T1.\\
$[X(\alpha,z) X(\delta_{\underline{m}},z),\  X(\beta,w)  X (\delta_{\underline{k}},w)]
= \begin{cases}0 \ \mbox{if}\  \alpha+\beta \ \notin \ \Delta \cup\{0\}\\
 N_{\alpha.\beta} \ X(\alpha + \beta,w)  X(\delta_{\underline{m}+{\underline{k}}},w) \ \delta(z/w)\ \\ \mbox{if}\ \alpha +\beta \ \in \ \Delta\\
\sigma_{\alpha} H(\alpha,w)  X(\delta_{\underline{m}+{\underline{k}}},w)  w  \delta(z/w)\\
-(X_{\alpha},X_{\beta})K \big(X(\delta_{\underline{m}+{\underline{k}}},w)  w  D_{z}  \delta(z/w)\\
- w  \delta_{\underline{m}} (w) X(\delta_{\underline{m}+{\underline{k}}},w)  \delta(z/w) \\
- X(\delta_{\underline{m}+{\underline{k}}},w) \ \delta(z/w)\big) \ \mbox{if}\ \alpha +\beta =0 \end{cases}$
$$
\begin{array}{llll}
T2.& [H(i,z) X(\delta_{\underline{m}},z),\  X(\alpha,w) X(\delta_{\underline{k}},w)]\\
&= \alpha (H(i)) \ X(\alpha,w)  X(\delta_{\underline{m}+{\underline{k}}},w)  w^{-1} \ \delta(z/w)\\[2mm]
T3.& [H(i,z) X(\delta_{\underline{m}},z), \ H(j,w) X(\delta_{\underline{k}},w)]\\
&= -(H(i), H(j)) K \big(X(\delta_{\underline{m}+{\underline{k}}},w)  w^{-1}  D_{z}  \delta(z/w)  -w^{-1}  \delta_{\underline{m}} (w)  X(\delta_{\underline{m}+{\underline{k}}},w)  \delta(z/w) \big)\\[2mm]
\end{array}
$$
We also need to check that the derivations act correctly and that will be verified at the end.
By definition it follows that $X(\delta_{\underline{m}},z), \delta_{\underline{k}} (z) X (\delta_{\underline{m}},z)$
commute with $\tau$ action.\\
First recall that $X_{k} (\alpha) = X(\alpha) \otimes t^{k}_{n}$ and $H_{k}(i) = H(i) \otimes t^{k}_{n} \ \in \gg_{aff}$
and they satisfy the following relations.\\

\noindent
A(1).\\
$\ [X_{m}(\alpha),\ X_{k} (\beta)]
=\begin{cases} 0 \ \mbox{if}\  \alpha+\beta \ \notin \ \Delta \ \cup \{0\}\\
 N_{\alpha.\beta} \ X_{m+k}(\alpha + \beta) , \alpha +\beta \ \in \ \Delta\\
\sigma_{\alpha} H_{m+k}(\alpha) +(X_{\alpha},X_{\beta}) m \ \delta_{m+k,0} \ K_{n} \ \mbox{if} \ \alpha+\beta =0\end{cases}$\\[3mm]
$A(2).\ [H_{m}(i),X_{k}(\alpha)]= \alpha \big(H(i)\big) \ X_{k+m} (\alpha)$\\[2mm]
$A(3).\ [H_{m}(i),H_{k}(j)]= \big(H(i),H(j)\big) \ m \ \delta_{m+k,0} \ K_{n}$\\

For $\alpha \in \Delta$ and $i=1,2,\cdots l$\\
Define\\
$X(\alpha,z) = \sum X_{k} (\alpha) \ z^{-k}$ and \\
$H(\alpha,z) = \sum H_{k} (\alpha) \ z^{-k-1}$.\\
Then the relations $A(1)$ to $A(3)$ can be written in the following infinite series.
\noindent
AS(1).\\
$[X(\alpha,z),\ X (\beta,w)]
= \begin{cases}0 \ \mbox{if}\  \alpha+\beta \ \notin \ \Delta \ \cup\{0\}\\
 N_{\alpha.\beta} \ X(\alpha + \beta,w) \ \mbox{if}\ \alpha +\beta \ \in \ \Delta\\
\sigma_{\alpha} H(\alpha,w) w \ \delta (z/w) - (X(\alpha),X(\beta))  K_{n}  z  D_{z}  \delta(z/w) \\ \ \mbox{if} \ \alpha+\beta =0\end{cases}$\\[2mm]
\noindent
$AS(2).\ [H(i,z),X(\alpha,w)]= \alpha \big(H(i)\big) \ X(\alpha,w) \ w^{-1}  \delta(z/w)$\\[2mm]
$AS(3).\ [H(i,z),H(j,w)]= -\big(H(i),H(j)\big) K_{n} w^{-1}  D_{z} \delta(z/w)$\\[3mm]
We will now check $T1,T2$ and $T3$.\\
Suppose $\alpha,\beta \in \Delta , \underline{m},\underline{k}\ \in \ZZ^{n-1}$. \\
Let $Y =[X(\alpha,z)  X(\delta_{\underline{m}},z), X(\beta,w)  X(\delta_{\underline{k}},w)]$.\\
Then $Y =[X(\alpha,z),  X(\beta,w)]  X(\delta_{\underline{m}},z)  X(\delta_{\underline{k}},w)$\\
If $\alpha+\beta \ \notin \ \Delta \ \cup\{0\}$ then it is clear $Y=0$\\
Suppose $\alpha + \beta \in \Delta$\\
Then $Y= N_{\alpha,\beta}  X(\alpha+\beta,w)  \delta(z/w). X(\delta_{\underline{m}},z)  X(\delta_{\underline{k}},w)$ (by AS1)\\
$=N_{\alpha,\beta}  X(\alpha+\beta,w)  X(\delta_{\underline{m}+\underline{k}},w)  \delta(z/w)$ \big(by Lemma (3.7)(1)\big)\\
This verifies second part of $T1$.\\[4mm]
Now suppose $\alpha+\beta=0$\\
Then 
$$
\begin{array}{llll}
Y&= \sigma_{\alpha} H(\alpha,w)\ w \delta(z/w)  X(\delta_{\underline{m}},z) X(\delta_{\underline{k}},w)\\
&-(X_{\alpha},X_{\beta}) K z D_{z}  \delta(z/w). X(\delta_{\underline{m}},z) X(\delta_{\underline{k}},w)\ \mbox{(By AS1)} \\
&= \sigma_{\alpha} H(\alpha,w) w \delta(z/w)  X(\delta_{\underline{m}},w) X(\delta_{\underline{k}},w)\\
&-(X_{\alpha},X_{\beta}) K [w  X(\delta_{\underline{m}},w) X(\delta_{\underline{k}},w) D_{z}  \delta(z/w)\\
&-w \ \delta_{\underline{m}}  (w)  X(\delta_{\underline{m}},w) X(\delta_{\underline{k}},w)  \delta(z/w)\\
&-X  (\delta_{\underline{m}} , w)  X(\delta_{\underline{k}},w) \delta(z/w)] \ \mbox{(by Lemma (3.7) and (3.9))}\\
&= \sigma_{\alpha} H(\alpha,w)\ w  X(\delta_{\underline{m}+{\underline{k}}},w) \delta(z/w)\\
&-(X_{\alpha},X_{\beta}) K [w  X(\delta_{\underline{m}+{\underline{k}}},w) D_{z}  \delta(z/w)\\
&-w  \delta_{\underline{m}}  (w)  X(\delta_{\underline{m}+{\underline{k}}},w) \delta(z/w)\\
&-X  (\delta_{\underline{m}+{\underline{k}}}, w)  \delta(z/w)]\\
\end{array}
$$
This proves the third part of $T1$.\\[4mm]
We will now verify $T2$.\\
Consider
$$
\begin{array}{llll}
T2.& [H(i,z) X(\delta_{\underline{m}},z), X(\alpha,w) X(\delta_{\underline{k}},w)]\\
&= [H(i,z), X(\alpha,w)] X(\delta_{\underline{m}},z)  X(\delta_{\underline{k}},w) \\
&= \alpha (H(i))  X(\alpha,w)  w^{-1} \ \delta(z/w)  X(\delta_{\underline{m}},z) \ X(\delta_{\underline{k}},w) \ \mbox{(By AS2)}\\
&= \alpha (H(i))  X(\alpha,w) \ w^{-1}  X(\delta_{\underline{m}+{\underline{k}}},w)\  \delta(z/w) \ \mbox{(By Lemma 3.7)}\\
\end{array}
$$
Which verifies $T2$.\\
We will now verify $T3$.\\
Consider
$$
\begin{array}{llll}
[H(i,z) X(\delta_{\underline{m}},z), H(j,w) X(\delta_{\underline{k}},w)]\\
= [H(i,z), H(j,w)]  X(\delta_{\underline{m}},z)  X(\delta_{\underline{k}},w)\\
= -(H(i), H(j)) K w^{-1}  D_{z}  \delta(z/w)  X(\delta_{\underline{m}},z)  X(\delta_{\underline{k}},w) \ \mbox{(By AS3)}\\
= -(H(i), H(j)) K w^{-1}  X(\delta_{\underline{m}},w)  X(\delta_{\underline{k}},w)  D_{z}  \delta(z/w)\\
+ (H(i), H(j)) K w^{-1}  \delta_{m} (w)  X(\delta_{\underline{m}},w)  X(\delta_{\underline{k}},w) \delta(z/w) \ \mbox{(By Lemma 3.9 (2))}\\
= -(H(i), H(j)) K [w^{-1}  X(\delta_{\underline{m}+{\underline{k}}},w) D_{z}  \delta(z/w)  - w^{-1}  \delta_{\underline{m}}  (w)  X(\delta_{\underline{m}+{\underline{k}}},w) \delta(z/w)].\\
\end{array}
$$
It is easy to verify $(1 \leq i \leq n-1)$
$$
[d_{i}, X \otimes t^{\overline{m}}] = m_{i} \ X \otimes t^{\overline{m}}, X \in \ \gg
$$
For the n th derivation consider
$$
\begin{array}{llll}
&[ {\overline{d}}_{n} \otimes 1 + 1 \otimes L_{0}, \displaystyle{\sum_{k}}\ X_{k} (\alpha)  X_{m_{n}-k} \ (\delta_{\underline{m}}) ]\\
&= \displaystyle{\sum_{k}} \ [{\overline{d}}_{n},X_{k} (\alpha)]  X_{m_{n}-k}  (\delta_{\underline{m}}) ]\\
&+ \displaystyle{\sum_{k}} \ X_{k} (\alpha) [L_{0}, X_{m_{n}-k}  (\delta_{\underline{m}})]\\
&= \displaystyle{\sum_{k}} \ k  X_{k} (\alpha) X_{m_{n}-k}  (\delta_{\underline{m}})\\
&+ \displaystyle{\sum_{k}} \ (m_{n} -k)  X_{k} (\alpha) X_{m_{n}-k}  (\delta_{\underline{m}})]\\
&= m_{n} \ \displaystyle{\sum_{k}} \  X_{k} (\alpha) X_{m_{n}-k}  (\delta_{\underline{m}}).\\
\end{array}
$$
\section{The module $V \otimes V (\Gamma)$}
In this section we analyize the module  $V \otimes V (\Gamma)$. Recall $V(\Gamma) = \CC [\Gamma] \otimes S (\aa_{-})$\\
Let $V (\lambda) = V \otimes \ e^{\lambda+Q} \otimes S (\aa_{-}), \lambda \in \Gamma$.\\
Which is clearly a $\stackrel{\sim}{\tau}$-module.\\
Further $V \otimes V (\Gamma) = \displaystyle{\bigoplus_{\lambda \in \Gamma/Q}} V (\lambda)$.\\
In this section we will prove that $V(\lambda)$ is isomorpic to $V(0)$ as $\stackrel{\sim}{\tau}$-module upto twist of an
automorphism. Consider $GL(n,\ZZ)$ the group of invertible matrices of order $n$ with entries in $\ZZ$. $GL(n,\ZZ)$ naturaly
act on $\ZZ^{n}$ and we denote this action as $B.\overline{m}$ for $\overline{m} \in \ZZ^{n}$ and $B \in GL(n,\ZZ)$\\
Let $\overline{m} = (m_{1}, \cdots, m_{n})$ and $B= (b_{ij})_{\stackrel{1 \leq i \leq n}{1 \leq j \leq n}}$\\
Let $B.\overline{m} = \overline{l}= (l_{1}, \cdots, l_{n}), l_{k} = \displaystyle{\sum_{i}} \ b_{ik} \ m_{i}$\\
We see that $GL(n,\ZZ)$ acts as an automorphisms on $\stackrel{\sim}{\tau}$ by 
$B.X \otimes t^{\overline{m}} = X \otimes t^{B.{\overline{m}}}, X \in \gg$\\
$B.d  (t^{\overline{m}}) t^{\overline{k}} = d (t^{B. \overline{m}})  t^{B. \overline{k}}$\\
$B.d_{i} = d_{i}^{1}, (d^{1}_{1},\cdots d^{1}_{n}) = B^{T^{-1}}(d_{1}, \cdots d_{n})$\\
Given $\lambda \in \Gamma$ define $B_{\lambda} = (b_{ij}) \in GL(n,\ZZ)$ by $b_{ii} = 1, b_{in}= (\lambda, \delta_{i}), i \neq n, b_{ij} =0$
otherwise. Let $e_{1},\cdots e_{n}$ be the standard basis of $\ZZ^{n}$. Then it is easy to check that 
$B_{\lambda} e_{i} = e_{i} + (\lambda, \delta_{i}) e_{n}, i \neq n, B_{\lambda}  e_{n} = e_{n}$. One can also check that\\
$$
\begin{array}{llll}
B_{\lambda} (X \otimes t^{\overline{m}})&=& X \otimes t^{\overline{m}}\  t^{(\lambda, \delta_{\underline{m}})}_{n}\\
B_{\lambda} (t^{\overline{m}} K_{i})&=& t^{\overline{m}} \ t_{n}^{(\lambda, \delta_{\underline{m}})}\ K_{i} + 
{(\lambda, \delta_{\underline{m}})} \ t^{\overline{m}} \ t_{n}^{(\lambda, \delta_{\underline{m}})}\ K_{n} \ \mbox{for}\ i \neq n\\
B_{\lambda} (t^{\overline{m}} K_{n})&=& t^{\overline{m}} \ t_{n}^{(\lambda, \delta_{\underline{m}})}\ K_{n}\\
\end{array}
$$
We now recall the following Lemma from $[EM]$ which is very standard and follows from the definition of vertex operators.
\paragraph{(5.1) Lemma} Let $\gamma \in \Gamma$ and $\delta \in Q, m \in \ZZ$\\
Then\\
\noindent
$X_{m} (\delta). e^{\gamma} \otimes 1 = \begin{cases} e^{\gamma+\delta} \otimes S_{-m-(\gamma,\delta)} (\delta) \ \mbox{if} \ m+(\gamma,\delta) < 0\\
e^{\gamma+\delta} \otimes 1 \ \mbox{if} \ m+(\gamma,\delta) = 0\\
0 \ \mbox{if} \ m+(\gamma,\delta) > 0 \end{cases}$\\
Where the operators $S_{p} (\delta)$ is defined by $exp \ T_{-} (\delta, z) = \displaystyle{\sum_{p=0}^{\infty}}  S_{p}  (\delta)  z^{p}.$ \qed \\
Note  that $S_{0} (\delta) = 1$.
Let $S_{p} (\delta) = 0, p< 0$.\\
Then $X_{m} (\delta)  e^{\lambda} \otimes 1 = e^{\lambda+\delta}  \otimes S_{-m-(\lambda,\delta)} (\delta)$\\
Denote the  $\stackrel{\sim}{\tau}$ module $V(\lambda)$ by $\pi_{\lambda}$.\\
Let $\varphi : V(\lambda) \rightarrow V(0)$ 
defined by $\varphi(v \otimes  e^{\lambda + \delta_{\underline{m}}} \otimes u) = v \otimes \ e^{\delta_{\underline{m}}} \otimes \ u$\\
which is a vector space isomorphism.
\paragraph{(5.2) Proposition :} $\pi_\gamma$ and $\pi_0$ are isomorphic as $\stackrel{\sim}{\tau}$ modules upto a twist of 
automorphism $B_\lambda$.\\
{\bf{Proof}} Consider for $X (\alpha) \in \gg, \overline{m} \in \ZZ^{n}, \delta \in Q, v \in V$ and 
$u \in S(\aa_{-}).$
$$
\begin{array}{llll}
(1)& \pi_{\lambda} \big(X(\alpha) \otimes t^{\overline{m}}\big) \ v \otimes e^{\lambda + \delta} \otimes u\\
&= \displaystyle{\sum_{k}} \ X_{k} (\alpha) \ v \otimes  X_{m_{n}-k} (\delta_{\underline{m}}) (e^{\lambda+\delta} \otimes u)\\
&= \displaystyle{\sum_{k}} \ X_{k} (\alpha) v \otimes \ e^{\lambda+\delta+\delta_{\underline{m}}} \ \otimes \ 
S_{-m_{n}+k-(\lambda,\delta_{\underline{m}})}(\delta_{\underline{m}}) \ u\\
\end{array}
$$
\noindent
\big( By Lemma 5.1 and the fact that $u$ commutes with $X_{m_{n}-k} (\delta) \big)$\\
Consider
$$
\begin{array}{llll}
(2)& \pi_{0} \ oB_{\lambda}. \ (X_{\alpha} \otimes t^{\overline{m}}) \ v \otimes e^{\delta} \otimes u\\
&= \pi_{0} \ o \ \big(X_{\alpha} \otimes t^{\overline{m}} \ t_{n}^{(\lambda,\delta_{\underline{m}})}\big) \ v \otimes e^{\delta} \otimes u\\
&= \displaystyle{\sum_{k}} \ X_{k} (\alpha) v \ \otimes \ X_{m_{n}-k+(\lambda,\delta_{\underline{m}})} (\delta_{\underline{m}}) (e^{\delta}  \otimes u)\\
&= \displaystyle{\sum_{k}} \ X_{k} (\alpha) v \ \otimes \ e^{\delta+\delta_{\underline{m}}} \otimes \  S_{-m_{n}+k-(\lambda,\delta_{\underline{m}})} \ (\delta_{\underline{m}}) \ u\\
\end{array}
$$
We see that (1) and (2) are equal upto the identification of $V(\lambda)$ and the $V(0)$ by $\varphi$.\\
The same calculation holds good for $H(i)$. Now we will check this for the center.\\
(3) Let $M= (\lambda,\delta_{\underline{m}})$ and $N= (\lambda,\delta_{\underline{k}})$\\
Consider
$$
\begin{array}{llll}
& \pi_{\lambda} \big(d (t^{\overline{m}}) t^{\overline{k}}\big) \ v \otimes e^{\lambda + \delta} \otimes u\\
&= \big( T_{m_{n}+k_{n}}^{\delta_{\underline{m}}} (\delta_{\underline{m}+{\underline{k}}}) +m_{n} X_{m_{n}+{k_{n}}} (\delta_{\underline{m}+{\underline{k}}}) \big)v \otimes \ e^{\lambda+\delta} \otimes u\\
&= v \otimes \displaystyle{\sum_{k}} \delta_{\underline{m}} (k)  X_{m_{n}+{k_{n}}-k} (\delta_{\underline{m}+{\underline{k}}}) \big) \ (e^{\lambda+\delta} \otimes u)\\
&+\  m_{n} \ v \otimes \ X_{m_{n}+{k_{n}}} (\delta_{\underline{m}+{\underline{k}}}) \big) \ (e^{\lambda+\delta} \otimes u)\\
&= v \otimes \displaystyle{\sum_{k \neq 0}} e^{\lambda+\delta+\delta_{\underline{m} + {\underline{k}}}} \otimes \delta_{\underline{m}} \ (k)
S_{-m_{n}-{k_{n}}-M-N+k} (\delta_{\underline{m}+{\underline{k}}}) u\\
&+\big(m_{n} + (\lambda,\delta_{\underline{m}})\big) v \otimes e^{\lambda+\delta+\delta_{\underline{m} + {\underline{k}}}} \otimes S_{-m_{n}-{k_{n}}-M-N} (\delta_{\underline{m}+{\underline{k}}}) u\\
\end{array}
$$
Now consider
$$
\begin{array}{llll}
(4)& \pi_{0} \ 0B_{\lambda}. \big(d (t^{\overline{m}}) t^{\overline{k}}\big) \ v \otimes e^{\delta}  \otimes u\\
&=\pi_{0} \big(d (t^{\overline{m}} \ t^{M}_{n}) t^{\overline{k}} \ t^{N}_{n}\big) \ v \otimes e^{\delta}  \otimes u\\
&= T_{m_{n}+{k_{n}+M+N}}^{\delta_{\underline{m}}} (\delta_{\underline{m}+{\underline{k}}}) (\ v \otimes e^{\delta}  \otimes u)\\
&+ \big(m_{n} + (\lambda,\delta_{\underline{m}})\big)  X_{m_{n}+{k_{n}+M+N}} (\delta_{\underline{m}+{\underline{k}}})  \ (v \otimes e^{\delta} \otimes u)\\
&= \displaystyle{\sum_{k}} \ v \ \otimes \delta_{\underline{m}} (k)   X_{m_{n}+{k_{n}+M+N-k}} (\delta_{\underline{m}+{\underline{k}}})  \ (e^{\delta} \otimes u)\\
&+ \big(m_{n} + (\lambda,\delta_{\underline{m}})\big)  X_{m_{n}+{k_{n}+M+N-k}} (\delta_{\underline{m}+{\underline{k}}})\ e^{\delta} \otimes u\\
&= \displaystyle{\sum_{k \neq 0}} \ v \ \otimes e^{\delta+\delta_{\underline{m}+{\underline{k}}}} \otimes  \delta_{\underline{m}} (k)
S_{-m_{n}-{k_{n}-M-N+k}} (\delta_{\underline{m}+{\underline{k}}})  \ u\\
&+ \big(m_{n} + (\lambda,\delta_{\underline{m}})\big) \ v \ \otimes e^{\delta+\delta_{\underline{m}}+{\underline{k}}} \otimes 
S_{-m_{n}-{k_{n}-M-N}} (\delta_{\underline{m}+{\underline{k}}})  \ u\\
\end{array}
$$
We omitted the term $k=0$ in the first term as $\delta_{\underline{m}} (0). e^{\delta} \otimes u = 0$. Now we see that 
(3) and (4) are equal upto the identification of $V(\lambda)$ and $V(0)$ via $\varphi$.\\
We will verify the action of the derivations.\\
Without loss of generality we can assume $\lambda = k_{1} d_{1} + \cdots + k_{n-1} d_{n-1} = d_{\underline{k}}$ so that 
$(\lambda,\lambda)=0$\\
Recall that the automorphism $B_{\lambda}$ act on $d_{i}$ as $(B^{T}_{\lambda})^{-1}$.\\
So that 
$$
\begin{array}{llll}
B_{\lambda}. d_{i} &= d_{i},  i \neq n\\
&=-(\lambda, \delta_{1})d_{1} - \cdots - (\lambda, \delta_{n-1}) d_{n-1} + d_{n} \ \mbox{if} \ i=n.
\end{array}
$$
It is easy to check that for $i \neq n$\\
$\pi_{\lambda} (d_{i}) = \pi_{0} \ B_{\lambda} (d_{i}) \ \varphi^{-1}$\\
Let $v$ in $V$ be such that $\overline{d}_{n} v =Pv$ for $P \in \CC$.\\
Let $v_{1} = v \otimes e^{\lambda + \delta_{\underline{m}}} \otimes a_{1} (-l_{1}) \cdots a_{d} (-l_{d})$\\
Let $v_{2} = v \otimes e^{\delta_{\underline{m}}} \otimes a_{1} (-l_{1}) \cdots a_{d} (-l_{d})$\\
Consider
$$
\begin{array}{llll}
\pi_{\lambda} d_{n}. v_{1} &=& (\overline{d}_{n} \otimes 1 + 1 \otimes L_{0}) v_{1}\\
&=& \overline{d}_{n} v \otimes e^{\lambda + \delta_{\underline{m}}} \otimes a_{1} (-l_{1}) \cdots a_{d} (-l_{d})\\
&&+ v \otimes L_{0} \big(e^{\lambda + \delta_{\underline{m}}} \otimes a_{1} (-l_{1}) \cdots a_{d} (-l_{d})\big)\\
&=& \big(P - (\lambda, \delta_{\underline{m}}) - \sum l_{i} \big) v_{1}.
\end{array}
$$
Now consider
$$
\begin{array}{llll}
\pi_{0} B_{\lambda} d_{n} v_{2} &= & \pi_{0} \big(-\displaystyle{\sum^{n-1}_{i=1}} (\lambda,\delta_{i}) d_{i} + \overline{d}_{n} \otimes 1 + 1 \otimes L_{0}\big) v_{2}\\
&=&\big(-(\lambda, \delta_{\underline{m}}) + P - \sum l_{i} \big) v_{2}
\end{array}
$$
This completes the verification.\\
Note that $v_{1}$ and $v_{2}$ are identified via the map $\varphi$.

\section{Integrability and finite dimensional weight spaces}
In this section we make remarks on integrability and on finite dimensional weight spaces.
\paragraph{(6.1) Remark} We would like to indicate that if $V$ is $\gg_{aff}$ integrable then $V \otimes V [\Gamma]$ is 
$\stackrel{\sim}{\tau}$- integrable in the Lie algebra case. In the case of super affine, most of the integrable modules are trivial 
and one need to consider partial integrable modules. We will not address this case here.\

So we suppose $\gg$ is simple finite dimensional Lie algebra. Recall that $\Delta$ is a root system of $\gg$ and 
$\stackrel{\sim}{\Delta} = \{ \alpha + \delta_{\underline{m}}, \delta_{\underline{k}}, \alpha \in \Delta \}$
is a root system of $\stackrel{\sim}{\tau}$. We call $\alpha + \delta_{\underline{m}}$ real if $\alpha \in \Delta$. We denote by 
$\stackrel{\sim}{\Delta}_{real}$ the set of all real roots. We call a $\stackrel{\sim}{\tau}$ module $W$ integrable if all real 
root vectors act locally nilpotently on $W$.\

We now take an integrable restricted module for $\gg_{aff}$ of non-zero level and with finite dimensional weight spaces. Now 
by Theorem (1.10) of $[E1]$ such a module is completely reducible. Thus we consider an irreducible integrable highest weight 
module $V$ for $\gg_{aff}$. (Highest weight follows because we are assumeing the module is restricted.) We will now indicate how 
$V \otimes  V (\Gamma)$ is integrable $\stackrel{\sim}{\tau}$-module. The argument are very standard and hence we will only 
sketch the proof.\

Consider the simple system $ \pi = \{\alpha_{1}, \cdots \alpha_{l}\}$ of $\gg$ and let $\beta$ be maximal root. Then let 
$\stackrel{\sim}{\pi} = \{ \alpha_{1},\alpha_{2} , \cdots \alpha_{l}, -\beta + \delta_{1}, \cdots -\beta + \delta_{n}\}$
which can be  thought of a simple system in the sense that $\ZZ-$linear span of $\stackrel{\sim}{\pi}$ is $\stackrel{\sim}{\Delta}$.
It is not too difficult to check that $ W. \stackrel{\sim}{\pi} = \stackrel{\sim}{\Delta}_{real}$  where the Weyl group $W$ is genarated
by reflections $\gamma_{\alpha_{1}}, \cdots \gamma_{\alpha_{l}}, \gamma_{-\beta+\delta_{1}}, \dots \gamma_{-\beta+\delta_{n}}.$
See $[E3]$ for more details. Recall that $V(\lambda) = V \otimes e^{\lambda+Q} \otimes S(\aa_{-})$ and in view of Proposition (5.2),
we can assume $\lambda = 0.$ It is easy to see that $v \otimes e^{0} \otimes 1$ genarates $V(\lambda)$ as		 
$\stackrel{\sim}{\tau}$- module where $v$ is the highest weight vector of $V$. Now by Lemma 3.4(b) of $[K1]$, to 
check nilpotency of a real root vector on $V(\lambda)$, it is sufficient to check on the genarator. Now by Proposition 3 of 
section 6.1 of $[MP]$ and the fact that $ W. \stackrel{\sim}{\pi} = \stackrel{\sim}{\Delta}_{real}$, it is sufficient to check
local nilpotency of operators $X_{0} (\alpha_{i}), i=1,2,\cdots l, X_{0} (-\beta +\delta_{j}) , j = 1,2,\cdots n-1$ and 
$X_{1} (-\beta)$ on the genarator $v \otimes e^{0} \otimes 1$. We already know that $X_{0} (\alpha_{i})$ is localy 
nilpotent as it acts on the first component.\\
Consider for $i \neq n$.\\
$X_{0} (-\beta + \delta_{i}) v \otimes e^{0} \otimes 1$\\
$=\displaystyle{\sum_{k}} \ X_{k} \ (-\beta)\  v \otimes  X_{-k} \ (\delta_{i}) \ (e^{0} \otimes 1)$\\
$=\displaystyle{\sum_{k \geq 0}} \ X_{k} \ (-\beta)\  v \otimes  X_{-k} \ (\delta_{i}) \ e^{0} \otimes 1
(\mbox{Since} \ X_{m} (\delta_{i}) . e^{0} \otimes 1 =0$ for $m>0$ \break by Lemma 5.1).\\
$=X_{0} (-\beta) v \otimes e^{\delta_{i}} \otimes 1$ (as $X_{k} (-\beta) v \otimes = 0$ for $k > 0$)\\
Now
$$
\begin{array}{llll}
(X_{0} (-\beta + \delta_{i})^{l} \ (v \otimes e^{0} \otimes 1) &= X_{0} (-\beta)^{l} \ v \otimes e^{\delta_{i}} \otimes 1\\
&= 0 \ \mbox{for} \ l > > 0.\\
\end{array}
$$
Now consider
$$
\begin{array}{llll}
X_{1} (-\beta)^{l} \ (v \otimes e^{0} \otimes 1) &=X_{1} (-\beta)^{l} v \  \otimes e^{0} \otimes 1\\
&= 0 \ \mbox{for} \ l > > 0.\\
\end{array}
$$
This completes the proof of integrability for the $\stackrel{\sim}{\tau}$- module $V \otimes V (\Gamma).$

\paragraph{(6.2) Remark} Let $V$ be a $\gg_{aff}$ module which is in the category $\mathcal{O}$. (See $[K1]$ for definition). 
Certainly $V$ is restricted. We assume that $V$ has finite dimensional weight space with respect to 
$\hh_{aff}= \hh \oplus \CC \ K_{n} \oplus \CC \ \overline{d}_n.$ Note that the eigenvalues of $\overline{d}_n$ are 
bounded above. We will now prove that $V \otimes V (\Gamma)$ has finite dimensional weight spaces with respect to 
$\stackrel{\sim}{\hh} = \hh \oplus \sum \ \CC K_{i} \oplus \ \sum \ \CC \ d_{i}$. Recall that the central element $K_{i}$ 
acts as $\delta_{i} (0)$ for $i\neq n$ and $K_{n}$ act as one. Consider $v \otimes e^{\lambda + \delta_{\underline{m}}} \otimes u, 
v \in V, u = a_{1} \ (-k_{1}) \cdots a_{d} \ (-k_{d}) \ \in S \ (\aa_{-})$ and assume that it is $\stackrel{\sim}{\hh}$ 
weight vector. We can assume $\lambda = l_{1} \ d_{1} + \cdots + l_{n-1} \ d_{n-1}$. Then by looking at the action of 
$K_{1}, \cdots K_{n-1}$ and $d_{1}, \cdots d_{n-1}$ we see that $\lambda$ and $\delta_{\underline{m}}$ are fixed. 
Let $v = \sum \ v_{i}, \ \overline{d}_n  v_{i} = k_{i} v_{i}$\\
Consider
$$
D=d_{n}. v_{i} \otimes \ e^{\lambda + \delta_{\underline{m}}} \otimes \ a_{1} \ (-l_{1}) \cdots a_{d} \ (-l_{d})
$$
and recall that $d_{n} = \overline{d}_{n} \otimes 1 + 1 \otimes \ L_{0}$ (See 3.5 for definition of $L_{0}$ and 
its action).\\
Then $D= N  \ v_{i} \otimes \ e^{\gamma + \delta_{\underline{m}}} \otimes \ a_{1} \ (-l_{1}) \cdots a_{d} \ (-l_{d})$\\
where $ N= k_{i} - \big( (\lambda,\delta_{\underline{m}}) + \sum l_{i} \big)$ which is fixed constant. As $k_{i}$ is bounded above 
and $-\displaystyle{\sum_{i}} \ l_{i}$ is bounded above we see that the possibilities for $k_{i}$ and $l_{i}$ are finite. 
This proves that the weight space is finite dimentional. In particular if $V$ is a $\hh_{aff}$ weight module than 
$V \otimes V (\Gamma)$ is a $\stackrel{\sim}{\hh}$ weight module. 

\paragraph{(6.3) Remark} Suppose $V$ and $W$ are restricted $\gg_{aff}$ - modules and suppose $f : V \rightarrow W$ a 
$\gg_{aff}$ - module map. Then there exists $\stackrel{\sim}{f} : V \otimes V (\Gamma) \longrightarrow W \otimes V (\Gamma)$ 
a $\stackrel{\sim}{\tau}$ - module map such that 
$\stackrel{\sim}{f} (v \otimes e^{\gamma} \otimes u) = f (v) \otimes e^{\gamma} \otimes u.$ 
Follows from the definition of $\stackrel{\sim}{\tau}$ - module.

\end{document}